\documentclass[xypic,syntonly,amssymb,verbatim,12pt]{amsart}
\usepackage{graphicx}
\usepackage[all]{xy}
\usepackage{amssymb}
\usepackage{amsmath}
\usepackage[mathscr]{euscript}
\usepackage{enumerate}

\usepackage{epsfig}
\theoremstyle{plain}
\newtheorem{Thm}{Theorem}
\newtheorem{Prop}{Proposition}

 \theoremstyle{definition}
\theoremstyle{remark}

\numberwithin{equation}{section}

\begin{document}
  \title {Gauge fields on $B$-branes over $\mathbb{CP}^n$}

 \author{ Andr\'{e}s   Vi\~{n}a} 
\address{Departamento de F\'{i}sica. Universidad de Oviedo.    Garc\'{\i}a Lorca 18.
     33007 Oviedo. Spain. }
\email{vina@uniovi.es}

 \maketitle
\begin{abstract}
Considering the $B$-branes over the complex projective space ${\mathbb P}^n$ as the objects of the bounded derived category $D^b({\mathbb P}^n)$, we  prove that the cardinal of the set of holomorphic gauge fields on a given $B$-brane
${\mathscr G}^{\bullet}$ is $\leq 1$. Moreover,   the cardinal is $1$ iff each ${\mathscr G}^p$ is 
isomorphic to a direct sum of copies of
 ${\mathscr O}_{{\mathbb P}^n}$.
\end{abstract}

\section {Introduction} \label{S:intro}

From the mathematical point of view, a $B$-brane over a complex manifold $X$ is an object of the bounded derived category $D^b(X)$ \cite{Aspin, Aspin-et}. 

 The definition of a holomorphic connection on a vector bundle proposed by Deligne \cite[page 6]{Deligne} can be generalized in a natural way to  objects of $D^b(X)$,  giving rise to the concept of a holomorphic field on $B$-branes \cite{Vina22}.  
In this context, the set ${\mathfrak c}({\mathscr F}^{\bullet})$ of gauge fields on the $B$-brane ${\mathscr F}^{\bullet}$, if it is nonempty, is an affine space whose vector space is 
${\rm Hom}_{D^b(X)}({\mathscr F}^{\bullet},\,\Omega^1_X{\mathscr F}^{\bullet})$ (see below, in Section \ref{S:setof gauge}).

 The purpose of this note is to study the cardinal  $\#\,{\mathfrak c}({\mathscr G}^{\bullet})$ of the set
$ {\mathfrak c}({\mathscr G}^{\bullet}),$  where ${\mathscr G}^{\bullet}$ is a
  $B$-brane over the complex projective space ${\mathbb P}^n$. In Theorem \ref{Th:card}, we prove that this cardinal is $\leq 1$, and 
  Theorem \ref{Th:2} gives a necessary and sufficient condition for this cardinal to be $1$.
  
  In the proofs of these  theorems plays an important role the fact that 
the set
\begin{equation}\label{set}
{\mathcal E}=\{{\mathscr  O}(-n),\, {\mathscr  O}(- n +1),\dots, {\mathscr O}(-1),\,{\mathscr O}\}
\end{equation}
is a full exceptional family in $D^b({\mathbb P}^n)$ \cite[Sect 8.3]{Huybrechts}. Thus, $D^b({\mathbb P}^n)$ is equivalent to the smallest triangulated subcategory that contains this exceptional family. This fact allows us to pass from  general branes to complexes whose terms are direct sums of elements    of the set
(\ref{set}).
 On these branes 
   the   functor
 ${\rm Hom}_{D^b({\mathbb P}^n)}\big(\,.\,,\,\Omega^1(\, .\,)\big)$ can be expressed in terms the cohomology of ${\mathbb P}^n$.
   This reduction permits us to prove the results stated in Theorems 
   \ref{Th:card} and \ref{Th:2}.  
   
Some well-known results in cohomology will be developed in some detail, in order to facilitate the understanding of the proofs of these theorems for those mathematical physicists who are not fully familiar with sheaf cohomology.


\section{The set of gauge fields}\label{S:setof gauge}

Given ${\mathscr G}^{\bullet}$ an object of $D^b({\mathbb P}^n)$, we set $\Omega^1({\mathscr G}^{\bullet}):=\Omega^1\otimes_{\mathscr O}{\mathscr G}^{\bullet}$, where $\Omega^1$ is the sheaf of holomorphic  $1$-forms on ${\mathbb P}^n$. By ${\mathscr J}({\mathscr G}^{\bullet})$ we denote the corresponding  first jet complex. That is, 
${\mathscr J}({\mathscr G}^{\bullet})$  is the  complex of {\it abelian sheaves}  $\Omega^1({\mathscr G}^{\bullet})\oplus {\mathscr G}^{\bullet}$ endowed with 
 the following ${\mathscr O}$-module structure
 \begin{equation}\label{O-ModEst}
 f\ast(\alpha^{\bullet}\oplus\sigma^{\bullet})=(f\alpha^{\bullet}+\partial f\otimes\sigma^{\bullet})\oplus f\sigma^{\bullet},
 \end{equation}
 where $f$ is an holomorphic function, $\alpha^{\bullet}\in {\Omega^1}({\mathscr G}^{\bullet})$ and $\sigma^{\bullet}\in {\mathscr G}^{\bullet}$.

  One has the following natural exact sequence of complexes of ${\mathscr O}$-modules, which in general {\it does not split} in the category of complexes of ${\mathscr O}$-modules,
  \begin{equation}\label{ExactSeq0}
  0\to \Omega^1({\mathscr G}^{\bullet})\to {\mathscr J}({\mathscr G}^{\bullet})\overset{\pi}{\rightarrow }{\mathscr G}^{\bullet}\to 0
  \end{equation}
  
 A {\it holomorphic gauge field} on the $B$-brane ${\mathscr G}^{\bullet}$  is a right inverse  in $D^b({\mathbb P}^n)$ of the morphism 
 $\pi$ \cite{Vina22}.  The existence of holomorphic gauge fields can be studied by considering certain 
 ${\rm Ext}$ groups.

 \smallskip
\noindent
{\sc Ext functors.} In general, given two bounded below  complexes $A^{\bullet}$ and  $B^{\bullet}$ in an abelian category ${\mathfrak C}$, the complex
${\rm Hom}^{\bullet}(A^{\bullet},\,B^{\bullet})$ is defined by (see \cite[page 17]{Iversen})
\begin{equation}\label{deltamg}
{\rm Hom}^m(A^{\bullet},\,B^{\bullet})=\prod_{p\in{\mathbb Z}}{\rm Hom}_{\mathfrak C}\big(A^p,\, B^{p+m}   \big),
\end{equation}
with the differential  $\delta_H$.
\begin{equation}\label{deltamg1}
(\delta_H^mg)^p=\delta_B^{m+p}g^p+(-1)^{m+1}g^{p+1}\delta_A^p.
\end{equation}

 If the category  ${\mathfrak C}$ has enough injectives, 
 $${\rm Ext}^i(A^{\bullet},\, B^{\bullet}):=H^i {\rm Hom}^{\bullet}(A^{\bullet},\,I^{\bullet}),$$
 where $I^{\bullet}$ is a complex of injective objects, which is $q$-isomorphic to $B^{\bullet}$. 
 
Some properties of funtors ${\rm Ext}^ i$  to be used here are the following:
 \begin{enumerate}[(a)]
 \item  ${\rm Ext}^0$ is the ${\rm Hom}$ functor of the   bounded below derived category $D^+(\mathfrak C)$. 
 \item If ${\mathfrak C}$ is the category of ${\mathscr O}$-modules and ${\mathscr G}^{\bullet}$ is a complex of locally free sheaves, then it is not necessary to take an injective resolution in the definition of functors 
 ${\rm Ext}^i({\mathscr G}^{\bullet},\, .\,)$. That is, ${\rm Ext}^i({\mathscr G}^{\bullet},\,{\mathscr F}^{\bullet})=H^i{\rm Hom}^{\bullet}({\mathscr G}^{\bullet},\,{\mathscr F}^{\bullet})$ \cite[Chap III, 6.5.1]{Hartshorne}.
 \item The family of functors ${\rm Ext}^i$  is  a $\delta$-functor \cite[page 205]{Hartshorne}.
\end{enumerate}
The fundamental properties of ${\rm Ext}$ functors can be found  \cite{Ge-Ma, Iversen, Kas-Sch}. 
 

 \smallskip
 
 \noindent
 {\sc The Atiyah class.}
 From (\ref{ExactSeq0}) together with the above property (c), it follows the Ext exact sequence
 \begin{align}\label{ExtExactSeq}
 0\to \,&{\rm Hom}_{D^b({\mathbb P}^n )}\big({\mathscr G}^{\bullet},\,\Omega^1( {\mathscr G}^{\bullet}  )    \big)\to
 {\rm Hom}_{D^b({\mathbb P}^n )}\big({\mathscr G}^{\bullet},\,{\mathscr J}( {\mathscr G}^{\bullet}  )    \big)\overset{\pi_*}{\longrightarrow} \\ \notag
\to \,&{\rm Hom}_{D^b({\mathbb P}^n )}\big({\mathscr G}^{\bullet},\, {\mathscr G}^{\bullet}   \big)\overset{\Delta}{\rightarrow}
 {\rm Ext}^1\big({\mathscr G}^{\bullet},\,\Omega^1( {\mathscr G}^{\bullet}  )    \big)\to\dots \notag
 \end{align}
 The element $\Delta({\rm id})\in {\rm Ext}^1\big({\mathscr G}^{\bullet},\,\Omega^1( {\mathscr G}^{\bullet}  )    \big)$ is the Atiyah class
 \cite{Atiyah}
 $a( {\mathscr G}^{\bullet})$ of ${\mathscr G}^{\bullet}$. 
 
 From the exactness of (\ref{ExtExactSeq}), we deduce:
 
 \begin{enumerate}
 \item
    The vanishing of $a( {\mathscr G}^{\bullet})$ is  a necessary and sufficient condition for the existence of a holomorphic gauge field on  the brane ${\mathscr G}^{\bullet}$.
    \item  The set of holomorphic gauge fields on this brane, if it is nonempty, is an affine space that has   ${\rm Hom}_{D^b({\mathbb P}^n )}\big({\mathscr G}^{\bullet},\,\Omega^1( {\mathscr G}^{\bullet}  )\big)$ as vector space.
    \end{enumerate}


 \subsection{Case ${\mathscr O}(k)$}\label{Ss:CaseSp}
 If ${\mathscr A}$ is an ${\mathscr O}$-module and $k$ is an integer, we use the habitual notation ${\mathscr A}(k)$ for the twisted ${\mathscr O}$-module 
${\mathscr A}\otimes_{\mathscr O}{\mathscr O}(k)$. 
The functor $F_{(k)}:=-\otimes_{\mathscr O}{\mathscr O}(k)$ is an automorphism of the category of ${\mathscr O}$-modules, since $F_{(k)}$ and $F_{(-k)}$ are inverse.  Moreover, given two ${\mathscr O}$-modules ${\mathscr A}$, ${\mathscr B}$, then 
\begin{equation}\label{k,-k}
{\rm Hom}_{\mathscr O}({\mathscr A}(k),\,{\mathscr B})\simeq {\rm Hom}_{\mathscr O}({\mathscr A},\,{\mathscr B}(-k)).
\end{equation} 

 The brane ${\mathscr G}^{\bullet}$, such that ${\mathscr G}^0={\mathscr O}(k)$ and ${\mathscr G}^i=0$ for $i\ne 0$,  will be denoted by ${\mathscr O}(k)$. As ${\mathscr O}(k)$ is a   locally free ${\mathscr O}$-module, by the above property (b)  
  \begin{equation}\label{Exti}
  {\rm Ext}^i( {\mathscr O}(k),\,{\mathscr F}^{\bullet})=H^i{\rm Hom}_{\mathscr O}\big({\mathscr O}(k),\,{\mathscr F}^{\bullet}\big)
 =H^i {\rm Hom}_{\mathscr O}\big({\mathscr O},\,{\mathscr F}^{\bullet}(-k)\big),
 \end{equation}
 where (\ref{k,-k})   has been used.

 In particular, if the complex ${\mathscr F}^{\bullet}$ reduces to $\Omega^1(k)$, then the differential $\delta_H$ defined in (\ref{deltamg1}) is trivial and  from (\ref{Exti}) together with the property (a), one deduces
 $${\rm Hom}_{D^b({\mathbb P}^n) }\big( {\mathscr O}(k),\,\Omega^1(k))= {\rm Hom}_{\mathscr O}({\mathscr O},\,\Omega^1) =\Gamma({\mathbb P}^n,\,\Omega^1)=H^0({\mathbb P}^n,\,\Omega^1).$$
 From the above item (2), together with  $H^0({\mathbb P}^n,\,\Omega^1)=0$ \cite[page 4]{O-S-S}, it follows the following proposition about the cardinal of the set of the holomorphic gauge fields on ${\mathscr O}(k)$.
 \begin{Prop} \label{Prop:1}
 $$\#\{\text{holomorphic gauge fields on}\,\,{\mathscr O}(k)\}\in\{0,\,1\}.$$
\end{Prop}

 Let us assume that  $\psi\in {\rm Hom}_{D^b({\mathbb P}^n)}\big({\mathscr O}(k),\,{\mathscr J}( {\mathscr O}(k))\big)$  
  is a holomorphic gauge field on ${\mathscr O}(k)$.
  As 
  \begin{align}
  &{\rm Hom}_{D^b({\mathbb P}^n)}\big({\mathscr O}(k),\,{\mathscr J}({\mathscr O}(k))\big)=H^0{\rm Hom}^{\bullet} \big({\mathscr O}(k),\,{\mathscr J}({\mathscr O}(k))\big)=  \notag\\ 
  &{\rm Hom}^{0} \big({\mathscr O}(k),\,{\mathscr J}({\mathscr O}(k))\big)   = {\rm Hom}_{\mathscr O} \big({\mathscr O}(k),\,\Omega^1(k)\oplus{\mathscr O}(k)\big) \notag
  \end{align}
 and $\psi$ 
   is a right  inverse of  $\pi$ (see (\ref{ExactSeq0})) in $D^b({\mathbb P}^n)$, it can be written as $\psi=\nabla\oplus {\rm id}$, with $\nabla:{\mathscr O}(k)\to \Omega^1(k)$. The property $\psi(f\sigma)=f\ast\psi(\sigma)$   
implies that $\nabla(f\sigma)=df\sigma+f\nabla(\sigma)$, for $f\in{\mathscr O}$. Hence, $\nabla$ can be considered as a connection, in the usual sense, on the   line bundle $L$  associated to ${\mathscr O}(k)$.

 The connection form of $\nabla$ in a local frame $e$,
 defined on an open $U$, is the holomorphic $1$-form $A$ determined by $\nabla(e)=A\otimes e$.
The curvature $F$ of this connection is a $2$-form  defined over the whole ${\mathbb P}^n$, which on $U$ is equal to 
 $\partial A$. Hence, $F$ is a global holomorphic $2$-form. The first Chern class of $L$ is defined by the form $\frac{i}{2\pi}F$. As $H^0({\mathbb P}^n, \,\Omega^2)=0$ \cite[page 4]{O-S-S}, it turns out that the Chern class vanishes; that is, $k=0$. 
 
On the other hand,
for the case $k=0$, the morphism
$$\psi:{\mathscr O}\to\Omega^1\oplus {\mathscr O},\;\;\sigma\mapsto \partial\sigma\oplus\sigma$$
satisfies $\psi(f\sigma)=f\ast\psi(\sigma).$  It is a holomorphic gauge field on ${\mathscr O}$. From the above arguments together with Proposition \ref{Prop:1}, it follows the  proposition.
\begin{Prop}\label{T:1} The cardinal  $\#\,{\mathfrak c}({\mathscr O}(k))=\delta_{0k}.$
\end{Prop}

 
 \subsection{General case}
 
 As usual, we denote with ${\mathscr F}^{\bullet}[l]$, with $l\in{\mathbb Z}$, the complex ${\mathscr F}^{\bullet}$ shifted $l$ to the left. 
   Given  ${\mathscr A}, {\mathscr B}$ elements of the generating set (\ref{set}),
  let us consider morphisms $h$ between ${\mathscr A}':={\mathscr A}[l]$ and
  ${\mathscr B}':={\mathscr B}[l']$. We denote by ${\rm Cone}(h)={\mathscr A}'[1]\oplus {\mathscr B}'$ the mapping cone of $h$
   \cite[page 154]{Ge-Ma}. We define ${\mathcal E}^{(1)}$ the set obtained adding to ${\mathcal E}$ the elements of the form ${\rm Cone}(h)$. Hence, an element of ${\mathcal E}^{(1)}$ is a complex whose term at a position $p$ is either  $0$, or ${\mathscr O}(k)$, or a direct sum of ${\mathscr O}(k_1)\oplus {\mathscr O}(k_2)$, with $-n\leq k,k_1,k_2\leq 0$.
   
    Repeating the process with the elements of ${\mathcal E}^{(1)}$, one obtains ${\mathcal E}^{(2)}$, etc. The objects of the triangulated subcategory generated by the set ${\mathcal E}$ are elements which belong to some ${\mathcal E}^{(m)}$.
  
  Therefore, an object   of the triangulated subcategory of $D^b({\mathbb P}^n)$  generated by the set (\ref{set}) is a complex $({\mathscr G}^{\bullet},\,d^{\bullet})$, where ${\mathscr G}^{p}$ is a sheaf  of the form 
 \begin{equation}\label{bigoplus}
 {\mathscr G}^{p}=\bigoplus_{i\in S_p}{\mathscr O}(k_{pi}),
 \end{equation}
 with $-n\leq k_{pi}\leq 0$ and $i$ varying in a finite set $S_p$. (When $i$ ``runs over the empty set'',  the direct sum is taken to be $0$).

 Since ${\mathscr G}^{\bullet}$ is a complex of locally free ${\mathscr O}$-modules
 \begin{align} &{\rm Hom}_{D^b({\mathbb P}^n)}\big({\mathscr G}^{\bullet},\,\Omega^1{\mathscr G}^{\bullet}\big)=
 H^0{\rm Hom}^{\bullet}\big({\mathscr G}^{\bullet},\,\Omega^1{\mathscr G}^{\bullet}\big)= \notag \\ \notag
 &\{g\in{\rm Hom}^{0}\big({\mathscr G}^{\bullet},\,\Omega^1{\mathscr G}^{\bullet}\big)\,|\,\delta_H g=0\},
\end{align}
where $\delta_H$ is the operator defined in (\ref{deltamg1}).
 Hence, according to (\ref{deltamg}), it follows
 $${\rm Hom}_{D^b({\mathbb P}^n)}\big({\mathscr G}^{\bullet},\,\Omega^1{\mathscr G}^{\bullet}\big)\subset
  {\rm Hom}^{0}\big({\mathscr G}^{\bullet},\,\Omega^1{\mathscr G}^{\bullet}\big)=
  \prod_p{\rm Hom}_{\mathscr O}({\mathscr G}^p,\,\Omega^1{\mathscr G}^p).$$
   By the additivity of the functor ${\rm Hom}_{\mathscr O}(\,.\,,\,.\,)$, it follows
  \begin{equation}\label{proplus}
  {\rm Hom}_{D^b({\mathbb P}^n)}\big({\mathscr G}^{\bullet},\,\Omega^1{\mathscr G}^{\bullet}\big)\subset
  \prod_p\bigoplus_{i,j}{\rm Hom}_{\mathscr O}\big( {\mathscr O}(k_{pi}),\,\Omega^1(k_{pj})  \big), \notag
   \end{equation}

 The summand ${\rm Hom}_{\mathscr O}\big( {\mathscr O}(k_{pi}),\,\Omega^1(k_{pj})  \big)$ is equal to
 $${\rm Hom}_{\mathscr O}\big( {\mathscr O},\,\Omega^1(k_{pj}-k_{pi})  \big)=\Gamma({\mathbb P}^n,\,\Omega^1(k_{pj}-k_{pi}))=0,$$
since $H^0({\mathbb P}^n,\,\Omega^1(k))=0$, for any  $k$. Therefore,
  ${\rm Hom}_{D^b({\mathbb P}^n)}\big({\mathscr G}^{\bullet},\,\Omega^1{\mathscr G}^{\bullet}\big)=0.$
  
One has the following theorem, which generalizes Proposition \ref{Prop:1}.
 \begin{Thm}\label{Th:card} If ${\mathscr G}^{\bullet}$ is a $B$-brane over ${\mathbb P}^n$, then
 $\#\,\mathfrak{c}({\mathscr G}^{\bullet})\leq 1.$
 \end{Thm}

 
 Let $\psi$  be a holomorphic gauge field on the above $B$-brane ${\mathscr G}^{\bullet}$. Then
  \begin{align}
  &\psi\in{\rm Hom}_{D^b({\mathbb P}^n) }\big( {\mathscr G}^{\bullet},\,{\mathscr J}({\mathscr G}^{\bullet} ) \big)=
  H^0{\rm Hom}^{\bullet} \big({\mathscr G}^{\bullet},\, {\mathscr J}({\mathscr G}^{\bullet})  \big)\subset \notag \\  
 &\prod_{p}{\rm Hom}_{\mathscr O} \big( {\mathscr G}^p,\, (\Omega^1{\mathscr G}^p \oplus {\mathscr G}^p) \big). \notag 
 \end{align}
Thus, $\psi$ determies a family $\{\psi^p:  {\mathscr G}^p \to \Omega^1{\mathscr G}^p \oplus {\mathscr G}^p\} $ of morphisms of ${\mathscr O}$-modules.
As $\psi$ is a right inverse of $\pi$, $ \psi^p=\nabla^p\oplus {\rm id}^p$, where 
$\nabla^p:{\mathscr G}^p \to \Omega^1{\mathscr G}^p$. Furthermore, as in Section \ref{Ss:CaseSp},  $\nabla^p$ is a holomorphic connection on ${\mathscr G}^p$, for any $p$.

 The trace
  ${\rm tr}(F^p)$ of curvature  of $\nabla^p$ is a holomorphic $2$-form on ${\mathbb P}^n$; as $H^0({\mathbb P}^n,\,\Omega^2)=0$, that trace vanishes. Hence, the first Chern class of vector bundle associated to the locally free  sheaf ${\mathscr G}^p$ vanishes.
  
   On the other hand, the first Chern class of ${\mathscr G}^p$ is the sum  
   $$\sum_{i} c_1({\mathscr O}(k_{pi})).$$
    This class is 0 iff $k_{pi}=0$ for all $i$, since the $k_{pi}\leq 0$. 
 
 Therefore, the existence of a holomorphic gauge field on the brane  ${\mathscr G}^{\bullet}$ defined in (\ref{bigoplus}) implies that 
 ${\mathscr G}^{\bullet}$ is a sequence of direct sum of copies of ${\mathscr O}$
 \begin{equation}\label{bigoplusimp}
 \dots \to \bigoplus_{i\in S_p}{\mathscr O}\overset{d^p}{\longrightarrow} \bigoplus_{i\in S_{p+1}}{\mathscr O}\to\dots
 \end{equation}
   Since ${\rm Hom}_{\mathscr O}({\mathscr O},\, {\mathscr O})\simeq {\mathbb C}$, the map $d^p$ is given by a constant complex matrix. 
   
 On the other hand, given the brane (\ref{bigoplusimp}),  if the set of indices $S_p$ has $m_p$ elements, on  $\oplus_{1}^{m_p}{\mathscr O}$ we define the map $\varphi^p$,   
   $$\varphi^p(\sigma_1\oplus\dots\oplus\sigma_{m_p})= \partial\sigma_1\oplus\dots\oplus\partial\sigma_{m_p}.$$ 
   It is a holomorphic connection on $\oplus_{1}^{m_p}{\mathscr O}$. Moreover, the family $\{\varphi^p\}$ is compatible with the ``constant'' differentials $d^p$. Thus, this family defines  a holomorphic gauge field  on the brane defined by the complex
   (\ref{bigoplusimp}). 
    Therefore, we have the following theorem.
 \begin{Thm}\label{Th:2}
 A $B$-brane on ${\mathbb P}^n$ admits a holomorphic gauge field iff it is isomorphic to a brane 
  of the form (\ref{bigoplusimp}).
 \end{Thm}


\end{document}